\documentclass[12pt]{article}

\usepackage{a4,graphicx,amsmath,amsfonts,amssymb,bbm,subfigure}
\usepackage{arydshln}
\usepackage{color}
\usepackage{url}

\newcommand{\R}{\mathbbm{R}}

\newtheorem{definition}{Definition}
\newtheorem{theorem}{Theorem}

\newtheorem{lemma}{Lemma}
\newtheorem{remark}{Remark}

\begin{document}
\begin{center}
{\bf \Large Equivalent Systems for Differential Equations \\[1ex] 
with Polynomially Distributed Delay}

\vspace{5mm}

{\large Roland~Pulch}

{\small Institute of Mathematics and Computer Science, 
University of Greifswald, \\
Walther-Rathenau-Str.~47, 17489 Greifswald, Germany. \\
ORCID: 0000-0002-6503-3300 $\;$ Email: {\tt roland.pulch@uni-greifswald.de}}  
\bgroup
\renewcommand\thefootnote{\fnsymbol{footnote}}%
\renewcommand\thempfootnote{\fnsymbol{mpfootnote}}%
\footnotetext[0]{The Version of Record of this manuscript has been published and 
is freely available in \\ 
Mathematical and Computer Modelling of Dynamical Systems, 20 Jul 2026, \\
\url{https://www.tandfonline.com/doi/10.1080/13873954.2026.2702655}}
\egroup
\end{center}

%%%%%%%%%%%%%%%%%%%%%%%%%%%%%%%%%%%%%%%%%%%%%%%%%%%%%%%%%%%%%%%%%%%%%%%%%%%%%
%%%                         Abstract                                      %%%
%%%%%%%%%%%%%%%%%%%%%%%%%%%%%%%%%%%%%%%%%%%%%%%%%%%%%%%%%%%%%%%%%%%%%%%%%%%%%
\begin{center}
{\bf Abstract}

\begin{tabular}{p{13cm}}
We consider delay differential equations with a polynomially distributed delay. 
We derive an equivalent system of delay differential equations, which includes just two discrete delays. 
The stability of the equivalent system and its stationary solutions are investigated. 
Alternatively, a Gaussian quadrature generates a discretisation of the integral, which describes the distributed delay in the original delay differential equation. 
This technique yields an approximate differential equation with multiple discrete delays. 
We present results of numerical computations, where initial value problems of the differential equations are solved by numerical time integration. 
Therein, the two approaches are compared. 

\bigskip


\bigskip

\end{tabular}
\end{center}

%%%%%%%%%%%%%%%%%%%%%%%%%%%%%%%%%%%%%%%%%%%%%%%%%%%%%%%%%%%%%%%%%%%%%%%%%%%%%
%%%                        Introduction                                   %%%
%%%%%%%%%%%%%%%%%%%%%%%%%%%%%%%%%%%%%%%%%%%%%%%%%%%%%%%%%%%%%%%%%%%%%%%%%%%%%

\section{Introduction}

Mathematical modelling by delay differential equations (DDEs) is used 
in many applications such as epidemiology, population dynamics, traffic flow, 
and others; see~\cite{brauer,erneux,sipahi-etal,smith}. 
A DDE includes discrete delays and/or distributed delays. 
Numerical methods were developed for the solution of initial value problems (IVPs) of DDEs with discrete delays; 
see~\cite{bellen-zennaro,das,mitsui}. 
In the case of a distributed DDE with a gamma distribution or a similar distribution, it is well known that an equivalent system of ordinary differential equations (ODEs) can be arranged; see~\cite{cassidy,guglielmi-hairer,nevermann-gros,teslya-wolkowicz}. 
Thus, numerical methods for IVPs of ODEs are applicable; see~\cite{butcher}. 

We investigate a DDE including a distributed delay, where the weight function of the distribution is a polynomial. 
Such a polynomial is present in a uniform distribution and a beta distribution, for example. 
We construct an equivalent system of DDEs with just two discrete delays. 
Up to our knowledge, this contribution represents the first work introducing this equivalent system in the case of a distributed delay. 
It is shown that the solutions of the differential equations with polynomially distributed delay and the equivalent system coincide. 
We examine the stability of the alternative system, where we focus on stationary solutions. 
Furthermore, a scaling of the system is formulated, which keeps the magnitude of its solutions small, so that error amplification and bad condition are reduced. 
We also generalise the proposed approach to systems of DDEs with multiple polynomially distributed delays. 

Alternatively, the integral, which defines a distributed delay, can be discretised by a quadrature method; see~\cite{sargood-etal}.
This discretisation yields a DDE with multiple discrete delays, whose solution represents an approximation of the solution for the distributed DDE. 
We employ a Gaussian quadrature associated with the weight function; see~\cite{quarteroni,stoer-bulirsch}; in this technique. 

Finally, we demonstrate the results of numerical computations using two DDE problems with distributed delay as test examples. 
First, a nonlinear system, which represents an epidemiological model, is considered, where the two approaches, i.e., equivalent DDE system and approximate DDE system, are compared. 
Second, a singularly perturbed problem is solved using the equivalent DDE system. 
All IVPs are solved by a numerical method for DDEs with multiple discrete delays; see~\cite{shampine-thompson}.

%%%%%%%%%%%%%%%%%%%%%%%%%%%%%%%%%%%%%%%%%%%%%%%%%%%%%%%%%%%%%%%%%%%%%%%%%%%%%
%%%                     Problem Definition                                %%%
%%%%%%%%%%%%%%%%%%%%%%%%%%%%%%%%%%%%%%%%%%%%%%%%%%%%%%%%%%%%%%%%%%%%%%%%%%%%%

\section{Delay differential equations}
\label{sec:problem-def}
In this section, we introduce classes of DDEs and the problem under investigation. 

%%%%%%%%%%%%%%%%%%%%%%%%%%%%%%%%%%%%%%%%%%%%%%%%%%%%%%%%%%%%%%%%%%%%%%%%%%%%%%
\subsection{Types of delay differential equations}
We outline different classes of systems of DDEs and some examples of their applications.
\begin{itemize}
\item[i)] 
DDEs that include multiple discrete delays 
\begin{equation} \label{dde-system-discrete}
y'(t) = f \Big( y(t) , y(t-\tau_1) , y(t-\tau_2) , \ldots, y(t-\tau_m) \Big) 
\end{equation}
with state variables $y : \R \rightarrow \R^k$ and 
right-hand side $f : \R^{(k+1)n} \rightarrow \R^k$. 
The delays $\tau_1,\tau_2,\ldots,\tau_m$ are 
positive, pairwise different constants 
or functions depending on time and/or the state variables. 
Such DDEs are applied in population models or recruitment models~\cite{brauer} and traffic flow models~\cite[p.~57]{sipahi-etal}.
\item[ii)] 
DDEs that include a distributed delay
\begin{equation} \label{dde-system-distributed}
y'(t) = f \Big( y(t) , \int_I v(y(t-\tau)) g(\tau) \; \mathrm{d}\tau \Big) 
\end{equation}
with right-hand side $f : \R^k \times \R \rightarrow \R^k$ 
and a function $v : \R^k \rightarrow \R$. 
The weight function $g : I \rightarrow \R$ specifies a 
distribution of the delay, where $I \subseteq [0,\infty)$ denotes the domain of the integration. 
DDEs of this form are also used in population models~\cite{brauer} and traffic flow models~\cite[p.~79]{sipahi-etal}. 
Distributed delays were also considered in reaction-diffusion equations~\cite{sargood-etal}.
\item[iii)]
Neutral DDEs include a derivative of the state variables at a previous time in the right-hand side. 
For example, systems of the form 
$$ y'(t) = f \Big( y(t) , y(t-\tau) , y'(t-\tau) \Big) $$
with a single delay~$\tau$ are of neutral type. 
Such DDEs were considered in an application related to Maxwell's equation 
in~\cite{ruehli-etal}, for example. 
\end{itemize}
%There are also some more complicated classes of DDEs. 
%For example, the right-hand side of a differential equation may include a term $y(\lambda t)$ for $t \ge 0$ with a constant $\lambda \in (0,1)$.
In this article, DDEs of type~(i) and~(ii) are employed.

%%%%%%%%%%%%%%%%%%%%%%%%%%%%%%%%%%%%%%%%%%%%%%%%%%%%%%%%%%%%%%%%%%%%%%%%%%%%%%
\subsection{Differential equations with distributed delay}
Now, DDEs with distributed delay are considered more detailed. 
We illustrate the approach using a single autonomous DDE 
of the form 
\begin{equation} \label{dde}
y'(t) = f \Big( y(t) , \int_I y(t-\tau) g(\tau) \; \mathrm{d}\tau \Big)
\end{equation}
with a continuous function $f : \R \times \R \rightarrow \R$. 
The interval $I$ is either the unbounded domain $I=[0,\infty)$ 
or a compact domain $I = [a,b]$, $0 \le a < b < \infty$. 
The bounds $a=\tau_{\min}$ and $b=\tau_{\max}$ represent a 
minimum delay and a maximum delay, respectively. 
The weight function~$g : I \rightarrow \R$ has the properties of 
a probability density function:

\newpage

\begin{itemize}
    \item[i)] $g$~is measurable,
    \item[ii)] $g(\tau) \ge 0$ \quad for all~$\tau \in I$, 
    \item[iii)] $\displaystyle \int_I g(\tau) \; \mathrm{d}\tau = 1$ \quad (normalisation).
\end{itemize}
The condition~(iii) is imposed without loss of generality, since a normalising constant of the weight function~$g$ could be included in the right-hand side~$f$. 
It is only required that the integral of~$g$ is finite.
In the case of a bounded interval~$I$, 
a distributed delay is also called a bounded delay.

%%%%%%%%%%%%%%%%%%%%%%%%%%%%%%%%%%%%%%%%%%%%%%%%%%%%%%%%%%%%%%%%%%%%%%%%%%%%%
\begin{remark} {\em
We study the case of an autonomous DDE~\eqref{dde}. 
Nevertheless, the following approach can be applied to non-autonomous DDEs 
as well. 
Furthermore, generalisations to systems of DDEs are outlined 
in Section~\ref{sec:generalisation}}. 
\end{remark} 
%%%%%%%%%%%%%%%%%%%%%%%%%%%%%%%%%%%%%%%%%%%%%%%%%%%%%%%%%%%%%%%%%%%%%%%%%%%%%

In the case of a gamma distribution, 
the probability density function reads as follows
$$ g(\tau) = C \, x^p \mathrm{e}^{- \tau / \theta} 
\qquad \mbox{for} \;\; \tau \ge 0 $$
with parameters $\theta,p > 0$ and a normalising constant $C>0$. 
Let $p \ge 1$ be an integer. 
In this case, the distribution is also called an Erlang distribution. 
A linear chain technique yields an equivalent system of ODEs associated with the DDE~\eqref{dde}, where the solution~$y$ 
is one component in a solution of the ODEs; 
see~\cite{cassidy,nevermann-gros}. 
Consequently, numerical methods for ODEs can be used; 
see~\cite{butcher}.

We study the case of a polynomial~$g$ on a compact interval $[a,b]$. 
Thus, the weight function exhibits the general form
\begin{equation} \label{g-polynomial}
g(\tau) = \sum_{i=0}^n \alpha_i \tau^i  
\end{equation}
with $\alpha_n \neq 0$. 
An important class of probability distributions with polynomial density is the beta distribution. 
The probability density function 
\begin{equation} \label{density-beta} 
g(\tau) = C \, (\tau-a)^p \, (b-\tau)^q 
\qquad \mbox{for} \;\; \tau \in [a,b] 
\end{equation}
depends on the parameters $p,q \ge 0$, 
whereas $C>0$ is a normalising constant. 
The uniform distribution ($g$ constant) represents the 
special case of $p=q=0$. 

Solutions of IVPs, which include (multiple) discrete delays 
with maximum delay $\tau_{\max}$, are characterised as follows. 

%%%%%%%%%%%%%%%%%%%%%%%%%%%%%%%%%%%%%%%%%%%%%%%%%%%%%%%%%%%%%%%%%%%%%%%%%%%%%
\begin{definition} 
Let $\phi : [-\tau_{\max},0] \rightarrow \R$ be a continuous history function. 
A continuous function $y : [-\tau_{\max},\infty) \rightarrow \R$ is a 
solution of an IVP of the DDE, if it holds that
\begin{itemize}
    \item[i)] $y$ is continuously differentiable on $(0,\infty)$ 
    and $\lim_{t \rightarrow 0+} y'$ exists, 
    \item[ii)] $y$ satisfies the DDE for $t > 0$, 
    \item[iii)] $y(t) = \phi(t)$ for $t \in [-\tau_{\max},0]$.
\end{itemize}
\end{definition} 
%%%%%%%%%%%%%%%%%%%%%%%%%%%%%%%%%%%%%%%%%%%%%%%%%%%%%%%%%%%%%%%%%%%%%%%%%%%%%

We also impose these conditions on a solution of an IVP of 
the distributed DDE~\eqref{dde}.
The maximum delay is $\tau_{\max}=b$ for a distributed delay on $I=[a,b]$. 
The following lemma examines the differentiability of the solution at $t=0$.

%%%%%%%%%%%%%%%%%%%%%%%%%%%%%%%%%%%%%%%%%%%%%%%%%%%%%%%%%%%%%%%%%%%%%%%%%%%%%
\begin{lemma} \label{lemma:regularity}
Let $I=[a,b]$ and $y : [-b,\infty) \rightarrow \R$ be a continuous function.  
If $y$ is smooth for $t>0$ and satisfies the distributed DDEs~\eqref{dde} 
with continuous weight function~$g$ on~$[a,b]$, 
continuous function $f$ on $\R \times \R$,
while $y(t) = \phi(t)$ for $t \in [-b,0]$ with continuous history function~$\phi$, 
then the right-sided limit of the derivative exists and satisfies
\begin{equation} \label{right-limit} 
\lim_{t \rightarrow 0+} y'(t) = 
f \Big( \phi(0), \int_a^b \phi(-\tau) g(\tau) \; \mathrm{d} \tau \Big) . 
\end{equation}
\end{lemma} 
%%%%%%%%%%%%%%%%%%%%%%%%%%%%%%%%%%%%%%%%%%%%%%%%%%%%%%%%%%%%%%%%%%%%%%%%%%%%%

{\em Proof}:

The DDE~\eqref{dde} includes the integral 
$$ J(t) = \int_a^b y(t-\tau) g(\tau) \; \mathrm{d} \tau , $$
which represents a parameter-dependent integral with parameter~$t$.
Since $y$ and $g$ are continuous, 
the function $J$ is also continuous is~$t$. 
We calculate the limit
\begin{align*} 
\lim_{t \rightarrow 0+} y'(t) & = 
\lim_{t \rightarrow 0+} f \left( y(t) , J(t) \right) 
\, = \, f \Big( \lim_{t \rightarrow 0+} y(t) , \lim_{t \rightarrow 0+} J(t) \Big) \\[1ex]
& = f \left( y(0) , J(0) \right) \, = \, 
f \Big( y(0) , \int_a^b y(-\tau) g(\tau) \; \mathrm{d}\tau \Big) . 
\end{align*}
The initial values determine $y(0) = \phi(0)$ and $y(-\tau) = \phi(-\tau)$. 
\hfill $\Box$

It is well-known that replacing a discrete delay by a distributed delay induces some regularisation. 
However, the derivative $y'$ typically has a singularity
if the history function $\phi$ is not differentiable at $t=0$. 
Only in the special case of a smooth function~$\phi$ and that 
$\lim_{t \rightarrow 0-} \phi(t)$ coincides with~\eqref{right-limit}, 
the derivative $y'$ is continuous.

%%%%%%%%%%%%%%%%%%%%%%%%%%%%%%%%%%%%%%%%%%%%%%%%%%%%%%%%%%%%%%%%%%%%%%%%%%%%%
%%%                     Equivalent System                                 %%%
%%%%%%%%%%%%%%%%%%%%%%%%%%%%%%%%%%%%%%%%%%%%%%%%%%%%%%%%%%%%%%%%%%%%%%%%%%%%%

\section{Equivalent systems}
\label{sec:system}
We investigate the existence of equivalent systems for the DDE~\eqref{dde} 
in the case of a polynomially distributed delay with~\eqref{g-polynomial}.

%%%%%%%%%%%%%%%%%%%%%%%%%%%%%%%%%%%%%%%%%%%%%%%%%%%%%%%%%%%%%%%%%%%%%%%%%%%%%
\subsection{Preliminaries}

The construction of an equivalent system is based on the following property. 

%%%%%%%%%%%%%%%%%%%%%%%%%%%%%%%%%%%%%%%%%%%%%%%%%%%%%%%%%%%%%%%%%%%%%%%%%%%%%
\begin{theorem} \label{thm:auxiliary-fcn}
Let $y : [-b,\infty) \rightarrow \R$ be continuous. 
We define
\begin{equation} \label{auxiliary-fcn}
    x_i(t) = 
    \int_a^b y(t-\tau) \tau^i \; \mathrm{d}\tau 
\end{equation}
for $i \ge 0$ and $t \ge 0$. 
It follows that $x$ is continuously differentiable and
$$ x_i'(t) = 
y(t-a) a^i - y(t-b) b^i 
+ i x_{i-1}(t) $$
for $i \ge 0$ and $t > 0$, fixing $x_{-1} \equiv 0$.
\end{theorem} 
%%%%%%%%%%%%%%%%%%%%%%%%%%%%%%%%%%%%%%%%%%%%%%%%%%%%%%%%%%%%%%%%%%%%%%%%%%%%%

{\em Proof}:

% Let $a=\tau_{\min}$ and $b=\tau_{\max}$. 
Using the substitution $s = t - \tau$, 
the integral in~\eqref{auxiliary-fcn} changes into the form
$$ x_i(t) = \int_{t-b}^{t-a} y(s) \, (t-s)^i \; \mathrm{d}s $$
for $i \ge 0$. 
In the case of $i=0$, 
the fundamental theorem of calculus yields 
$$ x_i'(t) = y(t-a) - y(t-b) = y(t-a) a^0 - y(t-b) b^0 . $$
The binomial formula implies for $i \ge 1$
\begin{align*} x_i(t) & = 
\int_{t-b}^{t-a} y(s) \, 
\sum_{k=0}^i \binom{i}{k} t^k (-s)^{i-k} \; \mathrm{d}s 
\; = \; \sum_{k=0}^i \binom{i}{k} t^k \int_{t-b}^{t-a} y(s) \, (-s)^{i-k} \; \mathrm{d}s \\
& = \int_{t-b}^{t-a} y(s) \, (-s)^i \; \mathrm{d}s 
+ \sum_{k=1}^i \binom{i}{k} t^k \int_{t-b}^{t-a} y(s) \, (-s)^{i-k} \; \mathrm{d}s . 
\end{align*}
Now, the integrands are continuous and independent of~$t$. 
Thus, differentiation with respect to~$t$ can be done 
using the fundamental theorem of calculus. 
We obtain with the product rule of differentiation
\begin{align*} x_i'(t) & = 
y(t-a) (a-t)^i - y(t-b) (b-t)^i \\
& \qquad 
+ \sum_{k=1}^i \binom{i}{k} k t^{k-1} \int_{t-b}^{t-a} y(s) \, (-s)^{i-k} \; \mathrm{d}s \\
& \qquad 
+ \sum_{k=1}^i \binom{i}{k} t^k \big( y(t-a) (a-t)^{i-k} - y(t-b) (b-t)^{i-k} \big) \\
& = y(t-a) a^i - y(t-b) b^i 
+ \sum_{k=1}^i \binom{i}{k} k t^{k-1} \int_{t-b}^{t-a} y(s) \, (-s)^{i-k} \; \mathrm{d}s .
\end{align*}
We examine the sum further
\begin{align*}
& \sum_{k=1}^i \binom{i}{k} k t^{k-1} \int_{t-b}^{t-a} y(s) \, (-s)^{i-k} \; \mathrm{d}s \\
= & \; i \sum_{k=1}^i \frac{(i-1)!}{(k-1)!(i-k)!} \, t^{k-1} \int_{t-b}^{t-a} y(s) \, (-s)^{i-k} \; \mathrm{d}s \\
= & \; i \sum_{k=0}^{i-1} \frac{(i-1)!}{k!(i-1-k)!} \, t^k \int_{t-b}^{t-a} y(s) \, (-s)^{i-1-k} \; \mathrm{d}s \\
= & \; i \int_{t-b}^{t-a} y(s) \sum_{k=0}^{i-1} \binom{i-1}{k} t^k (-s)^{i-1-k} \; \mathrm{d}s \, = \, i \int_{t-b}^{t-a} y(s) \, (t-s)^{i-1} \; \mathrm{d}s \, = \, i x_{i-1}(t) . 
\end{align*}
Hence the formula of Theorem~\ref{thm:auxiliary-fcn} is shown. 
\hfill $\Box$

%%%%%%%%%%%%%%%%%%%%%%%%%%%%%%%%%%%%%%%%%%%%%%%%%%%%%%%%%%%%%%%%%%%%%%%%%%%%%
\begin{remark} {\em
If the function~$y$ is continuously differentiable everywhere, 
then a proof of Theorem~\ref{thm:auxiliary-fcn} is straightforward 
by interchanging differentiation and integration. 
Theorem~\ref{thm:auxiliary-fcn} demonstrates that just the continuity of~$y$ 
is sufficient for the differentiability of functions~\eqref{auxiliary-fcn}. }
\end{remark} 
%%%%%%%%%%%%%%%%%%%%%%%%%%%%%%%%%%%%%%%%%%%%%%%%%%%%%%%%%%%%%%%%%%%%%%%%%%%%%

%%%%%%%%%%%%%%%%%%%%%%%%%%%%%%%%%%%%%%%%%%%%%%%%%%%%%%%%%%%%%%%%%%%%%%%%%%%%%
\subsection{Derivation of systems}
In the case of a polynomial weight function~\eqref{g-polynomial}, 
the integral of the DDE~\eqref{dde} becomes 
\begin{equation} \label{polynomial-decomposition} 
\int_a^b y(t-\tau) \, g(\tau) \; \mathrm{d}\tau = 
\sum_{i=0}^n \alpha_i \int_a^b y(t-\tau) \, \tau^i \; \mathrm{d} \tau 
= \sum_{i=0}^n \alpha_i x_i(t) 
\end{equation}
including the functions~\eqref{auxiliary-fcn}. 
Using Theorem~\ref{thm:auxiliary-fcn}, 
we arrange the following system of DDEs. 

%%%%%%%%%%%%%%%%%%%%%%%%%%%%%%%%%%%%%%%%%%%%%%%%%%%%%%%%%%%%%%%%%%%%%%%%%%%%%
\begin{definition} 
Let~\eqref{dde} be a DDE including a polynomially distributed delay, 
see~\eqref{g-polynomial}. 
An equivalent system of DDEs with two discrete delays is 
\begin{subequations} \label{equivalent-system}
\begin{align} 
    y'(t) & = f \Big( y(t) , \sum_{i=0}^n \alpha_i x_i(t) \Big) \label{equiv1} \\
    x_0'(t) & = y(t-a) - y(t-b) \label{equiv2} \\[1.5ex]
    x_i'(t) & = y(t-a) a^i - y(t-b) b^i + i x_{i-1}(t)
    \quad \mbox{for} \;\; i=1,\ldots,n . \label{equiv3} 
\end{align}
\end{subequations}
\end{definition} 
%%%%%%%%%%%%%%%%%%%%%%%%%%%%%%%%%%%%%%%%%%%%%%%%%%%%%%%%%%%%%%%%%%%%%%%%%%%%%

Only a single discrete delay is present in the case of $a=0$. 
In an IVP of the distributed DDE~\eqref{dde}, 
the initial conditions $y(t) = \phi(t)$ for $t \in [-b,0]$ are predetermined. 
Concerning the variables $x_i$, initial conditions are required 
only at $t=0$. 
Due to~\eqref{auxiliary-fcn}, it follows that
\begin{equation} \label{initial-values}
x_i(0) = \int_a^b \phi(-\tau) \tau^i \; \mathrm{d}\tau 
\end{equation}
for $i=0,1,\ldots,n$.
Thus the initial values are computable. 
Often constant initial values $y(t) = y_0$ for $t \le 0$ 
are employed. 
In this case, we obtain
\begin{equation} \label{initial-values-constant} 
x_i(0) = y_0 \, \frac{b^{i+1}-a^{i+1}}{i+1} 
\end{equation}
for $i=0,1,\ldots,n$. 

Now, numerical methods for IVPs of DDEs with a finite number of 
discrete delays can be used; see~\cite{bellen-zennaro}. 
A discontinuity of $y'$ at $t=0$ is propagated to discontinuities of higher-order derivatives at later time points. 
However, state-of-the-art numerical integrators consider these singularities appropriately; cf.~\cite[p.~89]{bellen-zennaro}. 
Consequently, the convergence order of the methods is also preserved in the presence of multiple discrete delays.

The definition~\eqref{auxiliary-fcn} shows that a discontinuity in a derivative $y^{(p)}$ implies a discontinuity in the derivative $x^{(p+1)}$. 
Hence, the state variables~$y$ and the auxiliary functions~$x$ share the singularities at the same time points.  
The introduction of auxiliary functions does not induce additional singularities in the system~\eqref{equivalent-system}.

The following theorem shows the equivalence of the differential equations. 

%%%%%%%%%%%%%%%%%%%%%%%%%%%%%%%%%%%%%%%%%%%%%%%%%%%%%%%%%%%%%%%%%%%%%%%%%%%%%
\begin{theorem} \label{thm:equivalence}
Let $y$ be a solution of an IVP of the DDE~\eqref{dde}. 
It follows that $y$ and the additional functions~\eqref{auxiliary-fcn} 
satisfy the system~\eqref{equivalent-system}. 
Vice versa, let $y,x_0,\ldots,x_n$ be a solution of an IVP of the system~\eqref{equivalent-system} with initial values $\phi$ for $y$ 
and~\eqref{initial-values} for~$x_i$. 
Consequently, the part~$y$ solves the DDE~\eqref{dde} 
with initial values~$\phi$.
\end{theorem} 
%%%%%%%%%%%%%%%%%%%%%%%%%%%%%%%%%%%%%%%%%%%%%%%%%%%%%%%%%%%%%%%%%%%%%%%%%%%%%

Proof:

We consider a solution~$y$ of an IVP of the DDE~\eqref{dde}. 
Thus~$y$ satisfies the first equation in~\eqref{equivalent-system} 
due to the definition of the additional functions~\eqref{auxiliary-fcn}. 
Theorem~\ref{thm:auxiliary-fcn} guarantees that the 
functions~\eqref{auxiliary-fcn} satisfy 
Eqns.~\eqref{equiv2} and~\eqref{equiv3}. 

Alternatively, we examine a solution $y,x_0,\ldots,x_n$ of an IVP of the system~\eqref{equivalent-system}. 
We consider Eqns.~\eqref{equiv2} and~\eqref{equiv3} for 
the unknowns $x=(x_0,x_1,\ldots,x_n)^\top$ and the fixed function~$y$. 
Hence, \eqref{equiv2}, \eqref{equiv3} represents a system of ODEs 
for $x$, where the right-hand side is continuous 
and continuously differentiable with respect to~$x$.
Thus solutions of IVPs are unique. 
Using the initial values~\eqref{initial-values}, 
the solution of~\eqref{equiv2}, \eqref{equiv3} are the 
functions~\eqref{auxiliary-fcn} due to Theorem~\ref{thm:auxiliary-fcn}. 
The relation~\eqref{polynomial-decomposition} shows that 
Eq.~\eqref{equiv1} is equivalent to the DDE~\eqref{dde}. 
\hfill $\square$

%%%%%%%%%%%%%%%%%%%%%%%%%%%%%%%%%%%%%%%%%%%%%%%%%%%%%%%%%%%%%%%%%%%%%%%%%%%%%
\begin{remark} {\em
In Theorem~\ref{thm:equivalence}, it is not required that solutions 
of IVPs of the DDE~\eqref{dde} are unique.}
\end{remark} 
%%%%%%%%%%%%%%%%%%%%%%%%%%%%%%%%%%%%%%%%%%%%%%%%%%%%%%%%%%%%%%%%%%%%%%%%%%%%%

%%%%%%%%%%%%%%%%%%%%%%%%%%%%%%%%%%%%%%%%%%%%%%%%%%%%%%%%%%%%%%%%%%%%%%%%%%%%%
\subsection{Stability}
We investigate the stability of the equivalent system and 
its stationary solutions. 
In the system~\eqref{equivalent-system}, the parts~\eqref{equiv2},~\eqref{equiv3} exhibit the form
\begin{equation} \label{dde-part-x} 
x'(t) = A x(t) + u ( y(t-a) , y(t-b ) ) 
\end{equation}
with variables $x = (x_0,x_1,\ldots,x_n)^\top$, a matrix 
\begin{equation} \label{matrix-part-x} 
A = \begin{pmatrix} 
0 & & & & \\
1 & 0 & & & \\
& 2 & \ddots & & \\
& & \ddots & 0 & \\
& & & n & 0 \\
\end{pmatrix} ,
\end{equation}
and a linear function~$u$. 
The rank of the matrix~\eqref{matrix-part-x} is~$n$.
The matrix~\eqref{matrix-part-x} has the sole eigenvalue 
$\lambda = 0$ with algebraic multiplicity $n+1$. 
There is only a single eigenvector $e=(0,\ldots,0,1)^\top$. 
Thus, the geometric multiplicity of the eigenvalue is one. 
It follows that the linear ODE 
\begin{equation} \label{autonomous-system}
    \hat{x}'(t) = A \hat{x}(t)
\end{equation} 
is Lyapunov-stable for $n=0$ and unstable for $n \ge 1$, 
see~\cite[p.~376]{braun}. 
Hence, the linear term $Ax$ induces an instability in the 
system~\eqref{equivalent-system}. 
However, error terms grow polynomially and not exponentially 
in time, since there is no eigenvalue with a positive real part. 
In view of this instability, IVPs of a system~\eqref{equivalent-system} 
should be solved with high accuracy requirements to reduce the growth of errors. 

The stationary solutions of a DDE are of the form $y(t) = y^*$ for 
all~$t$ with a constant~$y^*$.
A constant $y^*$ represents a stationary solution of~\eqref{dde}, 
if and only if the algebraic equation
$$ f (y^*,y^*) = 0 $$
is satisfied. 
The associated stationary solutions of the 
additional functions~\eqref{auxiliary-fcn} read as follows 
\begin{equation} \label{stationary-x} 
x_i^* = \int_a^b y^* \tau^i \; \mathrm{d}\tau 
= y^* \, \frac{b^{i+1}-a^{i+1}}{i+1} 
\end{equation}
for $i = 0,1,\ldots,n$, 
as in the calculation of~\eqref{initial-values-constant}. 

If $y^*$ is an asymptotically stable stationary solution, 
then we obtain the convergence
$$ \lim_{t \rightarrow \infty} x_i(t) = \lim_{t \rightarrow \infty} 
\int_a^b y(t-\tau) \tau^i \; \mathrm{d}\tau =  
\int_a^b \lim_{t \rightarrow \infty} y(t-\tau) \tau^i \; \mathrm{d}\tau = 
x_i^* $$
for $i \ge 0$ 
provided that the initial values of~$y$ are sufficiently close to $y^*$ 
and the initial values~\eqref{initial-values} are used for $x_i$.
An interesting question is whether this stationary solution 
of~\eqref{equivalent-system} is asymptotically stable or Lyapunov-stable. 
%Thus small perturbations of the initial values are considered. 

%%%%%%%%%%%%%%%%%%%%%%%%%%%%%%%%%%%%%%%%%%%%%%%%%%%%%%%%%%%%%%%%%%%%%%%%%%%%%
\begin{theorem} \label{thm:stability}
We suppose the existence and uniqueness of solutions of IVPs for all $t \ge 0$. 
Let $y^*$ be a Lyapunov-stable stationary solution of a DDE~\eqref{dde}. 
Consequently, $(y^*,x_0^*,\ldots,x_n^*)^\top$ with $x_i^*$ from~\eqref{stationary-x} is a stationary solution of the equivalent system~\eqref{equivalent-system}, 
which is unstable for $n \ge 2$. 
\end{theorem} 
%%%%%%%%%%%%%%%%%%%%%%%%%%%%%%%%%%%%%%%%%%%%%%%%%%%%%%%%%%%%%%%%%%%%%%%%%%%%%

Proof:

We suppose $y^* = 0$ without loss of generality. 
It follows that $x_i^*=0$ for all~$i$. 
The Lyapunov-stability implies that for an $\varepsilon > 0$, there is a $\delta \in (0,\varepsilon]$ 
such that the solution of an IVP is bounded by 
$| y(t) | < \varepsilon$ for all $t > 0$, 
if the initial function satisfies $| \phi(t) | < \delta$ for all $t \le 0$.

Let $z=(y,x_0,\ldots,x_n)^\top$ be a solution of the equivalent system~\eqref{equivalent-system} 
with initial functions $z_0=(\phi,\psi_0,\ldots,\psi_n)^\top$. 
%We assume that the stationary solution $(y^*,x_0^*,\ldots,x_n^*)^\top$ is Lyapunov-stable. 
%Using $\varepsilon$ from above, there is a $\tilde{\delta} \in (0,\delta]$ 
%such that the solution of~\eqref{equivalent-system} satisfies
%$\|z(t)\|_{\infty} < \tilde{\varepsilon}$ for all $t > 0$, 
%provided that the initial values are in a $\tilde{\delta}$-neigh\-borhood of zero. 
Let $\tilde{\delta} \in (0,\delta]$. 
If $z_0$ is in a $\tilde{\delta}$-neigh\-bor\-hood of zero with respect to the maximum norm, 
then the part $\phi$ is in a $\delta$-neigh\-bor\-hood of zero. 
In the right-hand side of~\eqref{equiv2}, \eqref{equiv3}, 
we estimate the terms 
$$ \left| y(t-a) a^i - y(t-b) b^i \right| \le 
(a^n + b^n) \sup_{s \ge -b} |y(s)| \le (a^n + b^n) \varepsilon $$
for $i=0,1,\ldots,n$ and $t \ge 0$. 
Thus, it holds that $\| u \|_{\infty} \le (a^n+b^n)\varepsilon$ 
for the function~$u$ in~\eqref{dde-part-x}. 
We consider the two systems~\eqref{dde-part-x} and~\eqref{autonomous-system} with same initial conditions. 
The matrix~\eqref{matrix-part-x} is nilpotent with $A^n \neq 0$ and $A^{n+1} = 0$. 
Thus, its matrix exponential reads as follows
$$ \mathrm{e}^{tA} = \sum_{j=0}^n \frac{t^j}{j!} \, A^j . $$
All terms are non-negative for $t \ge 0$ in this sum. 
We consider solutions of the system~\eqref{autonomous-system} in the case of $n \ge 2$. 
In each $\tilde{\delta}$-neighborhood of zero, 
there are initial conditions such that $\| \hat{x}(t) \|_{\infty} \ge K t^n \ge K t^2$ 
for all $t \ge 1$, 
where the constant $K>0$ depends on the initial values. 
In contrast, it holds that 
$$ \| x(t) - \hat{x}(t) \|_{\infty} \le t 
\sup_{s \in [0,t]} \| u(y(s-a),y(s-b)) \|_{\infty} \le 
 (a^n+b^n)\varepsilon \, t  $$
for any $t \ge 0$.
We obtain a chain of inequalities 
$$ K t^2 \le \| \hat{x}(t) \|_{\infty} 
\le \| x(t) - \hat{x}(t) \|_{\infty} + \| x(t) \|_{\infty} 
\le (a^n+b^n)\varepsilon \, t + \| x(t) \|_{\infty} . $$
It follows that
$$ (K t^2 - (a^n+b^n) \varepsilon t)  \le \| x(t) \|_{\infty} $$
for all $t \ge 1$. 
Now, $\| x(t) \|_{\infty}$ is not bounded for increasing~$t$. 
Hence, the stationary solution zero is not Lyapunov-stable. 
\hfill $\square$

%%%%%%%%%%%%%%%%%%%%%%%%%%%%%%%%%%%%%%%%%%%%%%%%%%%%%%%%%%%%%%%%%%%%%%%%%%%%%
\begin{remark} {\em
The above proof of Theorem~\ref{thm:stability} does not work in the 
case of \mbox{$n=1$}, because a linear growth of a solution $\hat{x}$ 
of~\eqref{autonomous-system} could be compensated in a solution~$x$ 
of~\eqref{dde-part-x} by the term~$u$. 
However, it is very unlikely that the term~$u$ balances the 
term $Ax$ for all positive times. 
Thus, we expect that the stationary solution is unstable in the 
case~$n=1$ as well.}
\end{remark} 
%%%%%%%%%%%%%%%%%%%%%%%%%%%%%%%%%%%%%%%%%%%%%%%%%%%%%%%%%%%%%%%%%%%%%%%%%%%%%

%%%%%%%%%%%%%%%%%%%%%%%%%%%%%%%%%%%%%%%%%%%%%%%%%%%%%%%%%%%%%%%%%%%%%%%%%%%%%
\subsection{Scaling}
\label{sec:scaling}
Concerning the equivalent system~\eqref{equivalent-system}, 
the additional functions~\eqref{auxiliary-fcn} include powers~$\tau^i$. 
It follows that the functions~\eqref{auxiliary-fcn} exhibit huge magnitudes for large delays. 
This effect may be harmful due to the amplification of errors. 
However, the disadvantageous effect can be avoided by scaling 
in the original DDE~\eqref{dde}. 

We introduce the new time variable $\tilde{t} = t/b$. 
Since it holds that $b=\tau_{\max}$, all involved delays are 
smaller than or equal to one in the new variable. 
In particular, it holds that $0 < \frac{a}{b} < 1$.
% Let $c=a/b$. 
Now, the DDE~\eqref{dde} is transformed into
\begin{equation} \label{dde-scaled}
\tilde{y}'(\tilde{t}) = 
b \, f \Big( \tilde{y}(\tilde{t}) , 
\int_{a/b}^1 \tilde{y}(\tilde{t}-\tau) \tilde{g}(\tau) \; \mathrm{d}\tau \Big) . 
\end{equation}
The new weight function~$\tilde{g}$ is just the linear transformation of~$g$ 
from the interval $[a,b]$ to the interval $[\frac{a}{b},1]$, 
followed by a normalisation of~$\tilde{g}$. 
It holds that
\begin{equation} \label{tilde-g} 
\tilde{g}(\tau) = C \, 
g \big( \tfrac{a}{b} + \tfrac{1-a/b}{b-a} (\tau - a) \big) 
\end{equation}
with a normalizing constant~$C$.
The degree of a polynomial~$g$ is invariant in this transformation. 

The equivalent system is arranged for the scaled DDE~\eqref{dde-scaled}. 
We obtain
\begin{align} \label{equivalent-system-scaled}
    \begin{split}
    \tilde{y}'(\tilde{t}) & = 
    b \, f \Big( \tilde{y}(t) , \sum_{i=0}^n \tilde{\alpha}_i \tilde{x}_i(\tilde{t}) \Big) \\
    \tilde{x}_0'(\tilde{t}) & = \tilde{y}(\tilde{t}-\tfrac{a}{b}) - \tilde{y}(\tilde{t}-1) \\[1.5ex]
    \tilde{x}_i'(\tilde{t}) & = \tilde{y}(\tilde{t}-\tfrac{a}{b}) \left( \tfrac{a}{b} \right)^i - \tilde{y}(\tilde{t}-1) \quad \mbox{for} \;\; i=1,\ldots,n ,
    \end{split}
\end{align}
where the coefficients~$\tilde{\alpha}_i$ result from the weight function~$\tilde{g}$ in~\eqref{tilde-g}.

%%%%%%%%%%%%%%%%%%%%%%%%%%%%%%%%%%%%%%%%%%%%%%%%%%%%%%%%%%%%%%%%%%%%%%%%%%%%%
%%%                   Generalisation to systems                           %%%
%%%%%%%%%%%%%%%%%%%%%%%%%%%%%%%%%%%%%%%%%%%%%%%%%%%%%%%%%%%%%%%%%%%%%%%%%%%%%

\section{Generalisations to systems}
\label{sec:generalisation}
A scalar DDE is considered in Section~\ref{sec:system}.  
Now, we generalise the approach to systems of DDEs 
with polynomially distributed delays. 

%%%%%%%%%%%%%%%%%%%%%%%%%%%%%%%%%%%%%%%%%%%%%%%%%%%%%%%%%%%%%%%%%%%%%%%%%%%%%
\subsection{Systems with one distributed delay}
We examine the system 
\begin{equation} \label{dde-system} 
y'(t) = f \Big( y(t) , \int_a^b v(y(t-\tau)) g(\tau) \; \mathrm{d}\tau \Big)
\end{equation} 
with solution $y = (y_1,\ldots,y_k)^\top$ and a continuous right-hand side 
$f : \R^k \times \R \rightarrow \R^k$. 
The weight function~$g$ has the form~\eqref{g-polynomial} and 
the same properties as in Section~\ref{sec:problem-def}. 
In addition, a continuous function $v : \R^k \rightarrow \R$ is included. 
A single component of the solution can be reproduced by 
$v(y(t-\tau)) = y_i(t-\tau)$ with some $i \in \{1,\ldots,k\}$, 
which is nearly identical to the situation investigated 
in Section~\ref{sec:system}. 
Another example is $v(y(t-\tau)) = y_i(t-\tau)y_j(t-\tau)$ 
with $i,j \in \{ 1,\ldots,k \}$ and $i \neq j$. 
This case occurs in some epidemiological models; see~\cite[p.~119]{brauer}. 

An equivalent system can be set up straightforward for the 
system of DDEs~\eqref{dde-system}. 
Again the principle is Theorem~\ref{thm:auxiliary-fcn}, 
where the continuous function~$y$ just has to be replaced by 
the continuous function $v(y)$. 
We obtain the equivalent system
\begin{align} \label{equivalent-system-system}
\begin{split}
    y'(t) & = f \Big( y(t) , \sum_{i=0}^n \alpha_i x_i(t) \Big) \\
    x_0'(t) & = v(y(t-a)) - v(y(t-b)) \\[1.5ex]
    x_i'(t) & = v(y(t-a)) a^i - v(y(t-b)) b^i + i x_{i-1}(t)
    \quad \mbox{for} \;\; i=1,\ldots,n  
\end{split}
\end{align}
including the additional functions
$$ x_i(t) = \int_a^b v(y(t-\tau)) \tau^i \; \mathrm{d}\tau 
\qquad \mbox{for} \;\; i=0,1,\ldots,n . $$
If initial values are predetermined by $y(t) = \phi(t)$ for $t \in [-b,0]$ 
with a continuous function $\phi : [-b,0] \rightarrow \R^k$, 
then an IVP of the system~\eqref{equivalent-system-system} is completed by, 
cf.~\eqref{initial-values}, 
\begin{equation} \label{initial-values-system}
x_i(0) = \int_a^b v(\phi(-\tau)) \tau^i \; \mathrm{d}\tau . 
\end{equation}
Thus, an equivalent system with the same properties is available for the system~\eqref{dde-system}.

%%%%%%%%%%%%%%%%%%%%%%%%%%%%%%%%%%%%%%%%%%%%%%%%%%%%%%%%%%%%%%%%%%%%%%%%%%%%%
\subsection{Systems with multiple distributed delays}
In the case of $d$ distributed delays, we apply a system
\begin{equation} \label{dde-system-multiple} 
y'(t) = f \Big( y(t) , 
\int_{a_1}^{b_1} v_1(y(t-\tau)) g_1(\tau) \; \mathrm{d}\tau , \ldots , 
\int_{a_d}^{b_d} v_d(y(t-\tau)) g_d(\tau) \; \mathrm{d}\tau \Big)
\end{equation} 
with solution $y=(y_1,\ldots,y_k)^\top$ and a continuous 
right-hand side $f : \R^k \times \R^d \rightarrow \R^k$. 
The continuity of the functions $v_1 , \ldots v_d : \R^k \rightarrow \R$ is also assumed. 
Each weight function $g_{\ell} : [a_{\ell},b_{\ell}] \rightarrow \R$ 
($0 \le a_{\ell}<b_{\ell}<\infty$) is a polynomial 
$$ g_{\ell}(\tau) = \sum_{i=0}^{n_{\ell}} \alpha_{i,\ell} \, \tau^i $$
for $\ell = 1,\ldots,d$. 
Now, each distributed delay induces a set of additional functions
$$ x_{i,\ell}(t) = \int_{a_{\ell}}^{b_{\ell}} v_{\ell}(y(t-\tau)) \tau^i 
\; \mathrm{d}\tau \qquad \mbox{for} \;\; i=0,1,\ldots,n_{\ell} $$
and $\ell = 1,\ldots,d$. 
The equivalent system of~\eqref{dde-system-multiple} reads as follows
\begin{align} \label{equivalent-system-system2}
\begin{split}
    y'(t) & = f \Big( y(t) , 
    \sum_{i=0}^{n_1} \alpha_{i,1} x_{i,1}(t) , \ldots , 
    \sum_{i=0}^{n_d} \alpha_{i,d} x_{i,d}(t) \Big) \\
    x_{0,1}'(t) & = v_1(y(t-a_1)) - v_1(y(t-b_1)) \\[1.5ex]
    x_{i,1}'(t) & = v_1(y(t-a_1)) a_1^i - v_1(y(t-{b_1})) b_1^i 
    + i x_{i-1,1}(t)
    \quad \mbox{for} \;\; i=1,\ldots,n_1 \\[0.5ex]
    & \vdots \\
    x_{0.d}'(t) & = v_d(y(t-a_d)) - v_d(y(t-b_d)) \\[1.5ex]
    x_{i,d}'(t) & = v_d(y(t-a_d)) a_d^i - v_d(y(t-{b_d})) b_d^i 
    + i x_{i-1,d}(t)
    \quad \mbox{for} \;\; i=1,\ldots,n_d . 
\end{split}
\end{align}
Let $b_{\max} = \max\{ b_1,\ldots,b_d \}$.
Given initial values $y(t) = \phi(t)$ for $t \in [-b_{\max},0]$ 
with a continuous function $\phi : [-b_{\max},0] \rightarrow \R^k$, 
the additional initial values of the equivalent 
system~\eqref{equivalent-system-system2} are obtained as 
in~\eqref{initial-values-system}. 
Furthermore, a scaling as in Section~\ref{sec:scaling} has to be done 
using $\tilde{t} = t / b_{\max}$.

Now, the system~\eqref{equivalent-system-system2} exhibits the 
total dimension
$$ k_{\rm eq} = k + d + \sum_{\ell=1}^d n_{\ell} . $$
Up to $2d$ discrete delays $a_{\ell},b_{\ell}$ 
are included in this system. 

%%%%%%%%%%%%%%%%%%%%%%%%%%%%%%%%%%%%%%%%%%%%%%%%%%%%%%%%%%%%%%%%%%%%%%%%%%%%%
%%%                   Systems from Quadrature                             %%%
%%%%%%%%%%%%%%%%%%%%%%%%%%%%%%%%%%%%%%%%%%%%%%%%%%%%%%%%%%%%%%%%%%%%%%%%%%%%%

\section{Approximate DDEs from quadrature}
\label{sec:quadrature}
An approach to solve an IVP of a DDE with distributed delay consists of 
replacing the involved integral by a quadrature formula. 
The DDE~\eqref{dde} changes into
\begin{equation} \label{dde-quadrature}
y'(t) = f \Big( y(t) , \sum_{k=1}^m \omega_k y(t-\tau_k) \Big)
\end{equation}
including the (pairwise different) nodes 
$\{ \tau_1, \tau_2, \ldots , \tau_m \} \subset I$ 
and the real-valued weights $\{\omega_1,\omega_2,\ldots,\omega_m\}$ 
belonging to a quadrature formula. 
% The quadrature scheme is dedicated to a particular weight function~$g$.  
Thus, Eq.~\eqref{dde-quadrature} represents a DDE with $m$~discrete delays. 
In~\cite{sargood-etal}, the integral was discretised by the 
(composite) Simpson rule. 
Initial values for the DDE~\eqref{dde} with distributed delay yield 
initial conditions for the DDE~\eqref{dde-quadrature} with discrete delays, 
as the nodes are located inside the interval~$I$. 

Now, numerical methods for IVPs of DDEs with discrete delays can be applied to~\eqref{dde-quadrature}.
An exact solution of a DDE~\eqref{dde-quadrature} represents an 
approximation of a solution of a DDE~\eqref{dde} with distributed delay. 
However, computation work increases with the number of discrete delays. 
Hence, the number~$m$ should be kept small. 

Gaussian quadrature can be used for weighted integrals 
\begin{equation} \label{weighted-integral} 
J(f) = \int_I f(t) \, g(t) \; \mathrm{d}t , 
\end{equation}
where the weight function $g : I \rightarrow \R_0^+$ 
satisfies some assumptions, see~\cite[p.~171]{stoer-bulirsch}. 
There is an associated Gaussian quadrature scheme for each 
traditional probability distribution, 
where the weight function represents a probability density function. 
Table~\ref{tab:gaussian-quadrature} lists some of these quadrature schemes. 
Each Gaussian quadrature scheme features an optimality property 
concerning polynomials~$f$ in~\eqref{weighted-integral}. 
Given a fixed number~$m$ of nodes, the Gaussian quadrature exhibits 
a maximum possible degree of polynomial exactness~$2m-1$, i.e., 
$$ \int_I t^i g(t) \; \mathrm{d}t = \sum_{k=1}^m \omega_k \tau_k^i 
\qquad \mbox{for} \;\; i=0,1,\ldots,2m-1 . $$
The Gaussian quadrature is the unique quadrature rule featuring this optimal polynomial exactness. 
Consequently, a different choice of nodes yields a quadrature formula with a lower polynomial exactness. 

In the case of a bounded interval $[a,b]$ and a polynomial~$g$, 
the Gaussian quadrature is convergent for any continuous integrand, which can be shown by the theorem of Weierstrass on polynomial approximation. 
However, the convergence may be slow if the integrand is not smooth. 
There is an error formula for a Gaussian quadrature with $m$~nodes; 
see~\cite[p.~180]{stoer-bulirsch}. 
This formula indicates a small error provided that the integrand is $2m$-times continuously differentiable.
%Thus, Gaussian quadrature methods are appropriate for the 
%discretisation of distributed delays in DDEs. 

However, a discontinuity of the derivative $y'$ at $t=0$ propagates to discontinuities of higher-order derivatives at later discrete time points. 
The number of discontinuities increases with the number of discrete delays. 
Thus, a DDE~\eqref{dde-quadrature} resulting from a quadrature rule features many singularities of the solution in the case of large numbers of nodes. 
Nevertheless, a solution of an IVP of a system~\eqref{dde-quadrature} becomes more and more smooth for increasing time. 
Hence, the error bound of the Gaussian quadrature soon applies after some time. 

In contrast, a DDE system~\eqref{equivalent-system} includes just two discrete delays. 
These delays induce discontinuities within derivatives, which are also present in the solution of the original DDE~\eqref{dde} with distributed delay. 
Thus, no additional singularities are introduced, 
which represents an advantageous property of the equivalent system.

%%%%%%%%%%%%%%%%%%%%%%%%%%%%%%%%%%%%%%%%%%%%%%%%%%%%%%%%%%%%%%%%%%%%%%%%%%%%%
\begin{remark} {\em
An advantage of Gausssian quadrature is that the procedure works for bounded as well as unbounded intervals~$I$. 
In both cases, the nodes are located in the interior of~$I$ and all the weights are positive.
However, we only apply the case of bounded intervals.}
\end{remark} 
%%%%%%%%%%%%%%%%%%%%%%%%%%%%%%%%%%%%%%%%%%%%%%%%%%%%%%%%%%%%%%%%%%%%%%%%%%%%%

%%%%%%%%%%%%%%%%%%%%%%%%%%%%%%%%%%%%%%%%%%%%%%%%%%%%%%%%%%%%%%%%%%%%%%%%%%%%%
\begin{table} 
\caption{Gaussian quadrature rules for weight functions associated to traditional probability distributions.\label{tab:gaussian-quadrature} }
\begin{center}
    \begin{tabular}{lcl}
    distribution & range~$I$ & name \\ \hline
    uniform distribution & $[a,b]$ & Gauss-Legendre quadrature \\
    beta distribution & $[a,b]$ & Gauss-Jacobi quadratute \\
    exponential distribution & $[0,\infty)$ & Gauss-Laguerre quadrure \\
    gamma distribution & $[0,\infty)$ & generalised Gauss-Laguerre quadrure
    \end{tabular}
\end{center}
\end{table}
%%%%%%%%%%%%%%%%%%%%%%%%%%%%%%%%%%%%%%%%%%%%%%%%%%%%%%%%%%%%%%%%%%%%%%%%%%%%%

%%%%%%%%%%%%%%%%%%%%%%%%%%%%%%%%%%%%%%%%%%%%%%%%%%%%%%%%%%%%%%%%%%%%%%%%%%%%%
%%%                       Numerical Results                               %%%
%%%%%%%%%%%%%%%%%%%%%%%%%%%%%%%%%%%%%%%%%%%%%%%%%%%%%%%%%%%%%%%%%%%%%%%%%%%%%

\section{Numerical results} \nopagebreak
\label{sec:numerical-results}
We examine numerical computations of two test examples in this section. 
The calculations were performed on a 
FUJITSU Esprimo P920 Intel(R) Core(TM) i7-9700 CPU with 3.00 GHz (8 cores) 
and operating system Microsoft Windows~10. 
All computations were executed using the software package
{\sc Matlab}~\cite{matlab}. 

%%%%%%%%%%%%%%%%%%%%%%%%%%%%%%%%%%%%%%%%%%%%%%%%%%%%%%%%%%%%%%%%%%%%%%%%%%%%%
\subsection{SIR model}
Epidemiological models often consist of ODEs or DDEs, 
whose solutions include populations of susceptibles~$S$, infected~$I$, and recovered~$R$. 
In~\cite{taylor-carr}, an $SIR$ model was investigated that features a birth rate and a death rate, where the loss of a temporary immunity is described by a single discrete delay. 
In~\cite{pulch-ecmi22}, the discrete delay was replaced by a distributed delay and numerical solutions were computed using the approach of Section~\ref{sec:quadrature}. 
Now, we apply this $SIR$ model without birth and death as a test example. 

Let $y_1=S$, $y_2=I$, $y_3=R$ be the three unknown populations. 
The system of DDEs with distributed delay reads as follows
\begin{align} \label{sir}
\begin{split}
y_1'(t) & = - \sigma y_1(t) y_2(t) 
+ \theta \int_a^b y_2(t-\tau) g(\tau) \; \mathrm{d}\tau \\
y_2'(t) & = \sigma y_1(t) y_2(t) - \theta y_2(t) \\
y_3'(t) & = \theta y_2(t) 
- \theta \int_a^b y_2(t-\tau) g(\tau) \; \mathrm{d}\tau 
\end{split}
\end{align}
with infection rate $\sigma > 0$ and recovery rate $\theta > 0$. 
It follows that $y_1+y_2+y_3$ is constant in time. 
We assume the normalisation $y_1+y_2+y_3=1$. 
The delay appears only in the component~$y_2$. 
We apply a beta distribution with weight function~\eqref{density-beta} and parameters $p=q=2$. 
Thus the associated polynomial~\eqref{g-polynomial} has degree $n=4$. 

We set the other parameters to $\sigma = 0.1$ and $\theta = 0.05$. 
The constant initial values 
$y_1(t)=0.99$, $y_2(t)=0.01$, $y_3(t)=0$
are predetermined for $t \le 0$. 
IVPs are solved in the time interval $[0,1000]$. 
We examine two cases with different ranges of the beta distribution:
\begin{itemize}
    \item Case~(i): $[a,b]=[30,150]$, 
    \item Case~(ii): $[a,b]=[150,250]$.
\end{itemize}
The system~\eqref{sir} has the form~\eqref{dde-system} with~$v(y)=y_2$. 
To generate the equivalent system~\eqref{equivalent-system-system} 
with discrete delays, 
we use functions of type~\eqref{auxiliary-fcn}
\begin{equation} \label{sir-x} 
x_i(t) = \int_a^b y_2(t-\tau) \tau^i \; \mathrm{d}\tau 
\end{equation}
for $i=0,1,2,3,4$. 
The differential equations for~\eqref{sir-x} are like
in~\eqref{equiv2},~\eqref{equiv3}.
The initial values follow from~\eqref{initial-values-constant}.
In our numerical simulations, the system of DDEs~\eqref{sir} is scaled as 
specified in Section~\ref{sec:scaling}. 
The time intervals are $[0,\frac{20}{3}]$ and $[0,4]$, respectively, in the solution of the scaled systems. 
Yet, the time axis is rescaled to the original time 
in graphical illustrations. 

\begin{figure}
 \centering
    \subfigure[populations]{
    \includegraphics[scale=0.48]{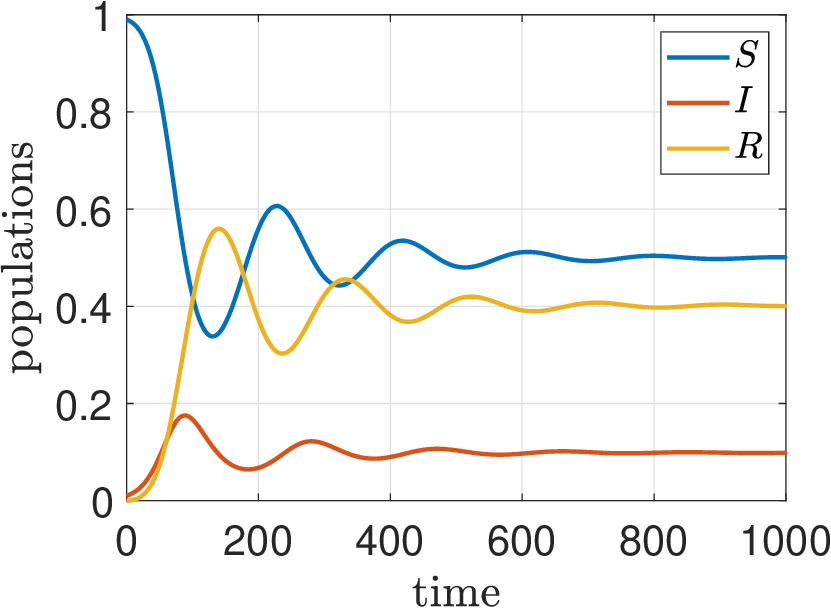}}
    \quad
    \subfigure[additional functions]{
    \includegraphics[scale=0.48]{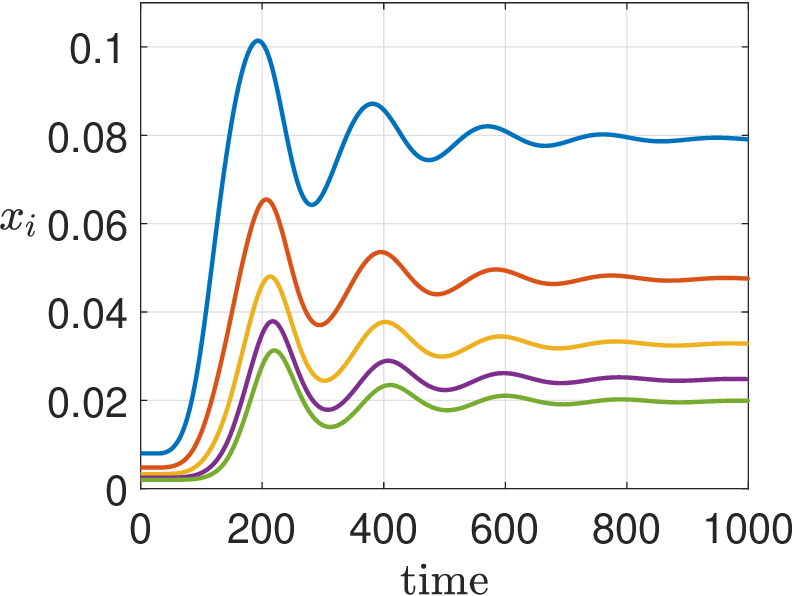}}
    \caption{Numerical solution of DDE system with two discrete delays in case~(i).}
    \label{fig:solution-a}
\end{figure}

\begin{figure}
 \centering
    \subfigure[populations]{
    \includegraphics[scale=0.48]{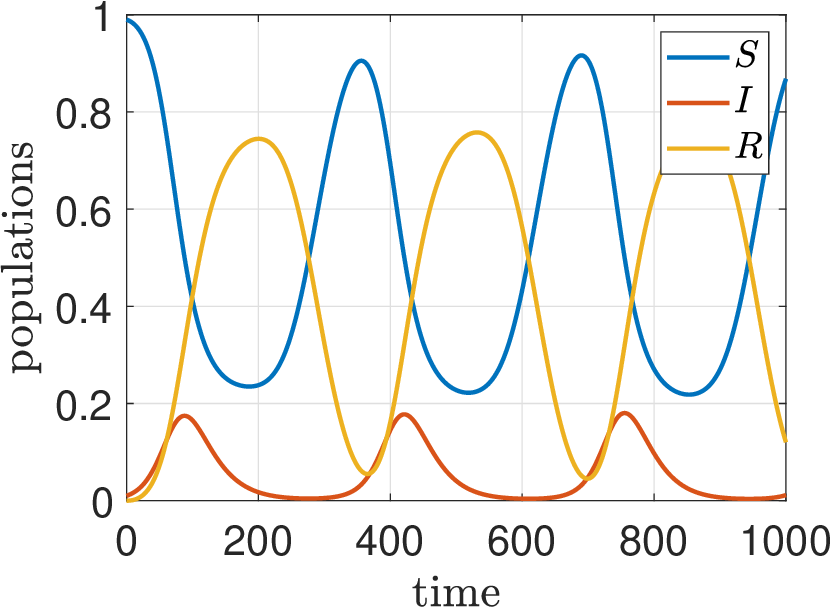}}
    \quad
    \subfigure[additional functions]{
    \includegraphics[scale=0.48]{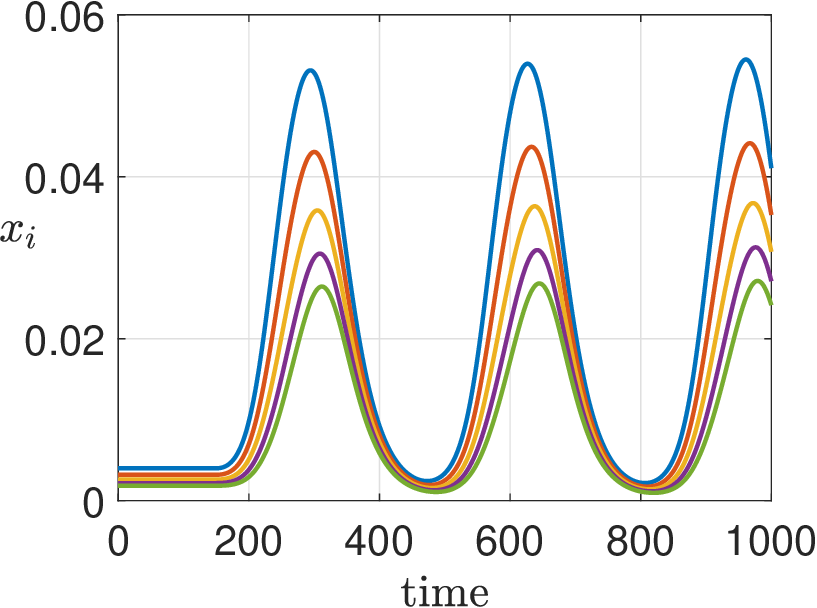}}
    \caption{Numerical solution of DDE system with two discrete delays in case~(ii).}
    \label{fig:solution-b}
\end{figure}

The built-in routine {\tt dde23}~\cite{shampine-thompson} of MATLAB, 
which implements an explicit method of Runge-Kutta type, 
produced the numerical solutions of IVPs of DDEs with discrete delays. 
This method features adaptive step size selection, 
where we require a relative error tolerance $\varepsilon_{\rm rel} = 10^{-6}$ 
and an absolute error tolerance $\varepsilon_{\rm abs} = 10^{-8}$. 
The integrator performs 610 steps in case~(i) and 1058 steps in case~(ii). 
Figure~\ref{fig:solution-a} and Figure~\ref{fig:solution-b} illustrate 
the solutions of the IVPs of the equivalent systems~\eqref{equivalent-system}
in the two cases. 
In case~(i), the solutions tend to a stationary solution. 
In case~(ii), the solutions tend to a periodic oscillation. 

%%% Figure: DDE solution %%%%%%%%%%%%%%%%%%%%%%%%%%%%%%%%%%%%%%%%%%%%%%%%%%%%
%%%%%%%%%%%%%%%%%%%%%%%%%%%%%%%%%%%%%%%%%%%%%%%%%%%%%%%%%%%%%%%%%%%%%%%%%%%%%

%%% Figure: DDE solution %%%%%%%%%%%%%%%%%%%%%%%%%%%%%%%%%%%%%%%%%%%%%%%%%%%%
%%%%%%%%%%%%%%%%%%%%%%%%%%%%%%%%%%%%%%%%%%%%%%%%%%%%%%%%%%%%%%%%%%%%%%%%%%%%%

\begin{figure}
 \centering
    \subfigure[case~(i)]{
    \includegraphics[scale=0.48]{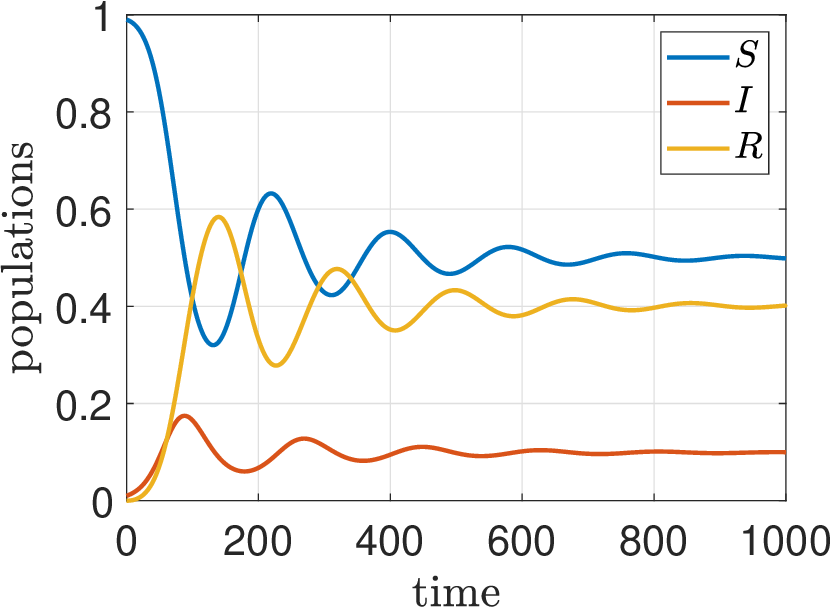}}
    \quad
    \subfigure[case~(ii)]{
    \includegraphics[scale=0.48]{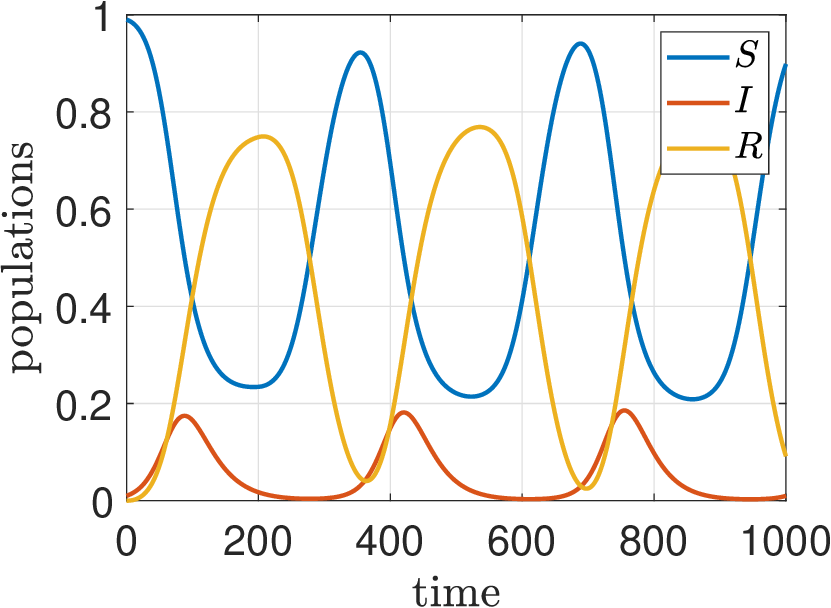}}
    \caption{Numerical solution of DDE systems including a single delay.}
    \label{fig:solution-quadrature}
\end{figure}

For comparison, we solve IVPs of the scaled approximate systems of 
DDEs~\eqref{dde-quadrature}, where Gaussian quadrature is used 
as described in Section~\ref{sec:quadrature}. 
The Gauss-Jacobi quadrature is associated to the beta distribution. 
In the case of a single node ($m=1$), we obtain a DDE with a single delay, 
which is identical to the expected value of the probability distribution. 
This instance represents an alternative model in comparison to 
the distributed delay. 
The expected values are $\bar{\tau}=90$ and $\bar{\tau}=200$ 
in case~(i) and case~(ii), respectively. 
The numerical solutions are depicted for the two cases 
in Figure~\ref{fig:solution-quadrature}. 
The solutions qualitatively agree to the situation of distributed delay. 
Only small quantitative differences are present. 

\begin{figure}
 \centering
    \subfigure[case (i)]{
    \includegraphics[scale=0.48]{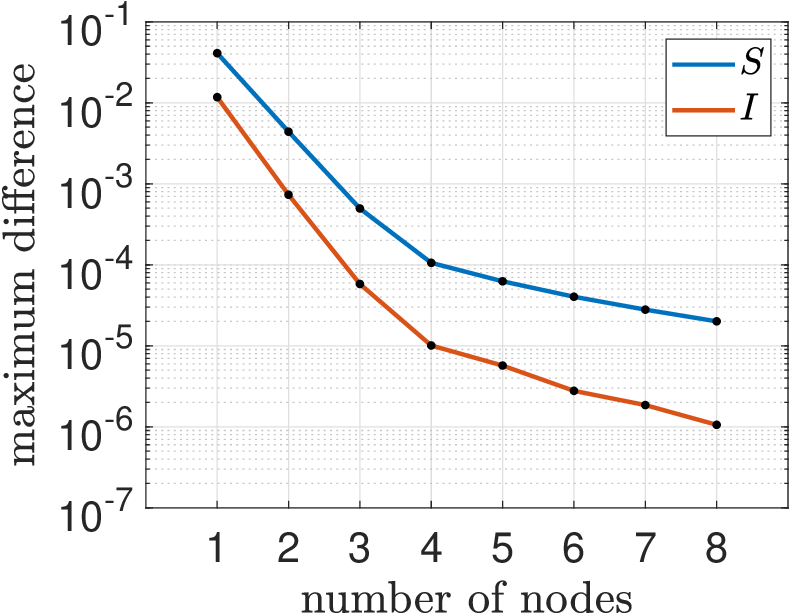}
    \label{fig:convergence-a}}
    \quad
    \subfigure[case (ii)]{
    \includegraphics[scale=0.48]{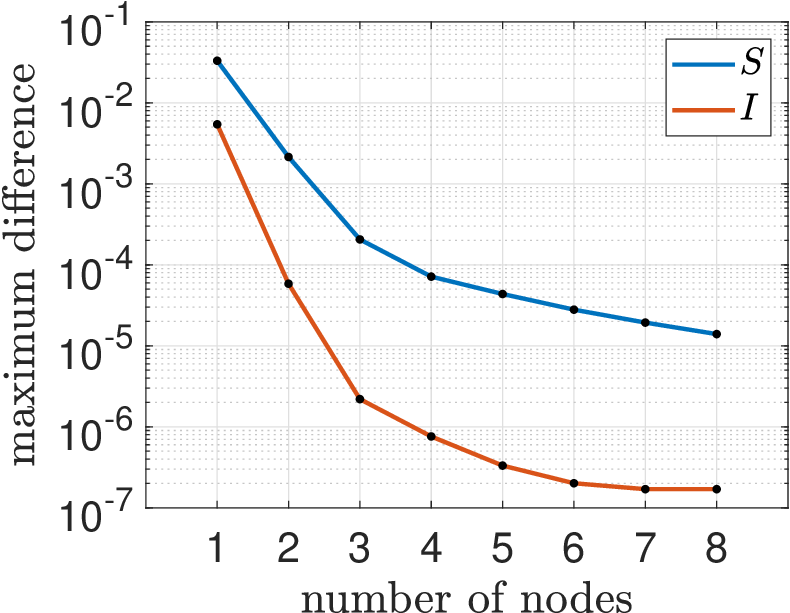}
    \label{fig:convergence-b}}
    \caption{Maximum difference between solutions of equivalent system and solutions of approximate systems using Gaussian quadrature with different numbers of nodes.}
    \label{fig:convergence}
\end{figure}

Finally, we determine the differences between the solution 
of the equivalent system and the solution 
of the approximate system using different 
numbers of nodes in the Gauss-Jacobi quadrature. 
In time integration, we impose a maximum step size 
$h_{\max} = 10^{-3}$ and $h_{\max} = 5 \cdot 10^{-4}$ 
in case~(i) and case~(ii), respectively. 
The solutions are evaluated at $1000$ equidistant points in 
the total time intervals 
$[0,\frac{20}{3}]$ and $[0,4]$, respectively. 
Figure~\ref{fig:convergence} shows the maximum differences for 
components~$S$ and~$I$, while the differences for $R$  
are almost identical to the differences for~$S$. 
We observe that the maximum differences exponentially decay
for increasing numbers of nodes. 
However, the differences stagnate for larger numbers of nodes, 
because the differences are dominated by the errors of the time integration 
if the quadrature errors are tiny. 
These results verify that the equivalent system of the form~\eqref{equivalent-system} 
yields the correct solution of the original DDE~\eqref{dde} in this example.

%%% Figure: DDE solution %%%%%%%%%%%%%%%%%%%%%%%%%%%%%%%%%%%%%%%%%%%%%%%%%%%%
%%%%%%%%%%%%%%%%%%%%%%%%%%%%%%%%%%%%%%%%%%%%%%%%%%%%%%%%%%%%%%%%%%%%%%%%%%%%%

%%% Figure: Convergence %%%%%%%%%%%%%%%%%%%%%%%%%%%%%%%%%%%%%%%%%%%%%%%%%%%%%
%%%%%%%%%%%%%%%%%%%%%%%%%%%%%%%%%%%%%%%%%%%%%%%%%%%%%%%%%%%%%%%%%%%%%%%%%%%%%

We also examine the computational effort in the time integration of the 
systems of DDEs. 
Table~\ref{tab:computation} depicts the number of integration steps 
as well as the computation times. 
We recognise that the number of steps, 
which is chosen by the local error control under the restriction of a maximum step size, 
is nearly the same in all systems used. 
As expected, the computation time grows for increasing numbers of discrete delays in the systems, since more delays cause more computation work. 
Consequently, the numerical integration of the equivalent transformed system of DDEs is cheaper than a numerical integration of a system applying Gaussian quadrature with a large number of nodes. 

%%%%%%%%%%%%%%%%%%%%%%%%%%%%%%%%%%%%%%%%%%%%%%%%%%%%%%%%%%%%%%%%%%%%%%%%%%%%%
\begin{table} 
\caption{Number of steps and computation times (seconds) in time integration of DDEs using Gaussian quadrature with different numbers of nodes  and transformed DDEs.\label{tab:computation}}
\begin{center}
    \begin{tabular}{ccccc}
    & case~(i) & & case~(ii) & \\ 
    no. nodes &  no. steps & comp. time & no. steps & comp. time \\ \hline
    1 & 6669 & 0.9 & 8002 & 1.2 \\
    2 & 6672 & 1.2 & 8006 & 1.5 \\
    3 & 6674 & 1.4 & 8004 & 1.8 \\
    4 & 6680 & 1.7 & 8010 & 2.2 \\
    5 & 6681 & 2.0 & 8017 & 2.5 \\
    6 & 6710 & 2.2 & 8047 & 2.9 \\
    7 & 6718 & 2.5 & 8062 & 3.3 \\
    8 & 6744 & 2.8 & 8076 & 3.6 \\ \hdashline
    transf. system & 6672 & 1.3 & 8006 & 1.7 
    \end{tabular}
\end{center}
\end{table}
%%%%%%%%%%%%%%%%%%%%%%%%%%%%%%%%%%%%%%%%%%%%%%%%%%%%%%%%%%%%%%%%%%%%%%%%%%%%%

%%%%%%%%%%%%%%%%%%%%%%%%%%%%%%%%%%%%%%%%%%%%%%%%%%%%%%%%%%%%%%%%%%%%%%%%%%%%%
\subsection{Singulary perturbed problem}
\begin{figure}
 \centering
    \includegraphics[scale=0.5]{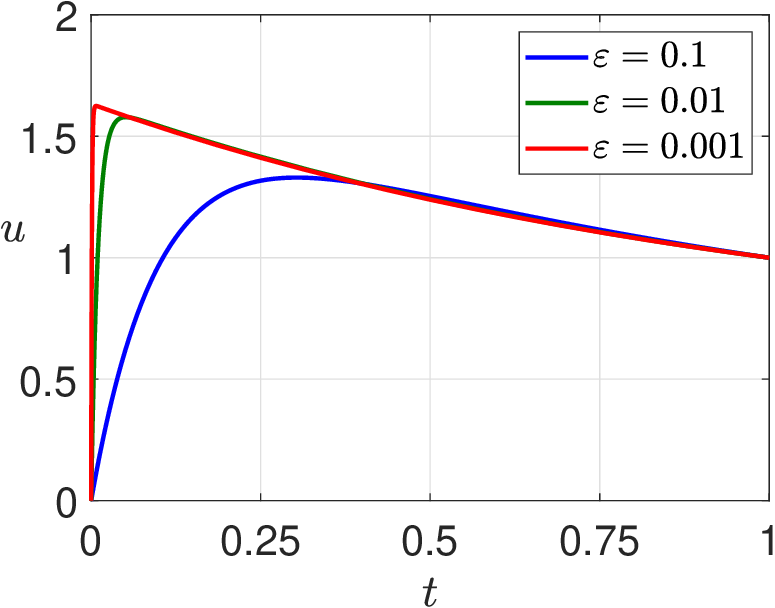}
    \caption{Solutions of singularly perturbed problem for different 
    parameters~$\varepsilon$ in the case of ODEs.}
    \label{fig:spp-ode}
\end{figure}

In the case of ODEs, 
a benchmark of a singulary perturbed problem is the 
boundary value problem (BVP) 
\begin{equation} \label{spp-ode}
\begin{array}{l}
\varepsilon u''(t) + u'(t) = - \mathrm{e}^{-t} 
\qquad \mbox{for} \;\; 0 < t < 1, \\[0.5ex]
u(0) = 0, \quad u(1) = 1, \\
\end{array}
\end{equation} 
with a parameter $\varepsilon > 0$. 
We solved this problem using the shooting method, 
where the unknown initial value $u'(0)$ was determined approximately by a bisection method. 
Therein, IVPs of the equivalent first-order system were solved by the 
MATLAB built-in routine {\tt ode23t}, 
which implements the (implicit) trapezoidal rule that includes 
a selection of step sizes based on a local error control. 
We set the error tolerances to 
$\varepsilon_{\rm rel} = 10^{-4}$ and $\varepsilon_{\rm abs} = 10^{-6}$. 
An implicit integration scheme was applied, 
since ODEs become stiff for small parameters~$\varepsilon$. 
Figure~\ref{fig:spp-ode} illustrates the solutions of the 
singularly perturbed problem~\eqref{spp-ode} for some parameters~$\varepsilon$. 
Furthermore, Table~\ref{tab:spp-ode} shows the initial values computed 
for the derivative by the shooting method as well as the number of 
integration steps to solve one IVP using this initial value. 
The results indicate the behaviour 
$u'(0) \rightarrow \infty$ for $\varepsilon \rightarrow 0$. 

%%% Figure: Singularly Perturbed Problem - ODE %%%%%%%%%%%%%%%%%%%%%%%%%%%%%%
%%%%%%%%%%%%%%%%%%%%%%%%%%%%%%%%%%%%%%%%%%%%%%%%%%%%%%%%%%%%%%%%%%%%%%%%%%%%%

%%%%%%%%%%%%%%%%%%%%%%%%%%%%%%%%%%%%%%%%%%%%%%%%%%%%%%%%%%%%%%%%%%%%%%%%%%%%%
\begin{table} 
\caption{Initial values of derivative and number of integration steps in singulary perturbed problem of ODEs.\label{tab:spp-ode}}
\begin{center}
    \begin{tabular}{crc}
    $\varepsilon$ & $u'(0)$ & no. steps \\ \hline
    $10^{-1}$ & 15.9 & 115 \\
    $10^{-2}$ & 162.9 & 157 \\
    $10^{-3}$ & 1631.9 & 197 \\
    $10^{-4}$ & 16323.2 & 231\\
    $10^{-5}$ & 163238.2 & 265 \\
    $10^{-6}$ & 1632457.9 & 298
    \end{tabular}
\end{center}
\end{table}
%%%%%%%%%%%%%%%%%%%%%%%%%%%%%%%%%%%%%%%%%%%%%%%%%%%%%%%%%%%%%%%%%%%%%%%%%%%%%

Now, we modify the singularly perturbed problem~\eqref{spp-ode} 
to a case of a differential equation featuring a distributed delay
\begin{equation} \label{spp-dde}
\begin{array}{l}
\varepsilon u''(t) + \tfrac{1}{b-a} \int_a^b u'(t-s) \; \mathrm{d}s 
= - \mathrm{e}^{-t} 
\qquad \mbox{for} \;\; 0 < t < 1, \\[0.5ex]
u(0) = 0, \quad u(1) = 1, \\
\end{array}
\end{equation} 
with a parameter $\varepsilon > 0$. 
The inherent delay exhibits a uniform distribution 
in the interval $[a,b]$. 
We employ the associated system of first order with $y_1 = u$ and $y_2 = u'$, 
i.e., 
\begin{align} \label{spp-dde-firstorder}
\begin{split}
y_1'(t) & = y_2(t) \\ 
y_2'(t) & = \frac{1}{\varepsilon} \left( 
- \frac{1}{b-a} \int_a^b y_2(t-s) \; \mathrm{d}s - \mathrm{e}^{-t} \right) 
\end{split}
\end{align}
with boundary values $y_1(0) = 0$ and $y_1(1) = 1$. 
Since DDEs are considered, initial values also have to be specified 
for $t < 0$ using a history function. 
In the following, we always fix constant initial values 
$y_i(t) = y_i(0)$ for each $t < 0$ and $i=1,2$ 
in the system~\eqref{spp-dde-firstorder}. 
Again, the BVP was solved by the shooting method, 
where a bisection determined the unknown initial value $y_2(0)$. 
However, a high accuracy is required in the bisection method 
to achieve a sufficiently accurate final value $y_1(1) \approx 1$. 
Thus, many iteration steps were performed to enclose appropriate initial values~$y_2(0)$. 

We applied the equivalent system~\eqref{equivalent-system} 
associated to~\eqref{spp-dde-firstorder}, 
which represents a DDE system for $y_1,y_2,x_0$ including 
the discrete delays $a$ and $b$. 
The interval $[a,b] = [0.2,0.4]$ was selected. 
Hence, scaling is not required in this example. 
We solved the equivalent system using the MATLAB routine {\tt dde15s}, 
which is a part of the supporting information to~\cite{agrawal-etal}. 
This routine implements an implicit multi-step method for systems of DDEs 
with discrete delays. 
%However, it turns out that the DDE system~\eqref{spp-dde-firstorder} 
%is less stiff than the equivalent ODE system of~\eqref{spp-ode} 
%for same parameter~$\varepsilon$. 
We fixed the error tolerances $\varepsilon_{\rm rel} = 10^{-6}$ and 
$\varepsilon_{\rm abs} = 10^{-8}$ in the selection of step sizes. 

\begin{figure}
 \centering
    \subfigure[solution~$y_1$]{
    \includegraphics[scale=0.48]{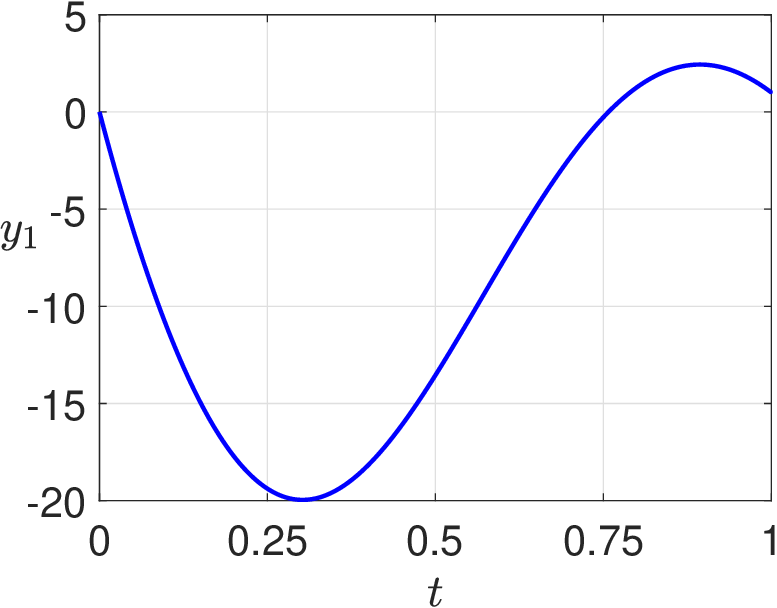}}
    \quad
    \subfigure[solution~$y_2$]{
    \includegraphics[scale=0.48]{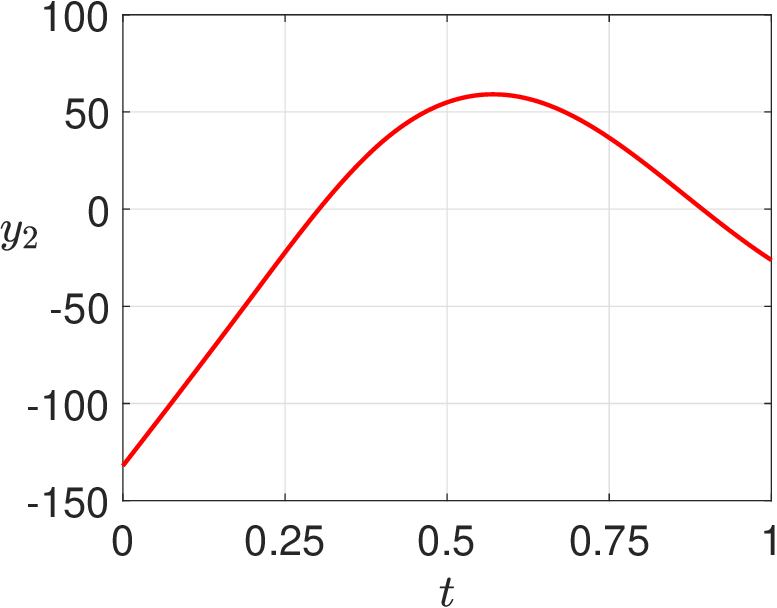}}
    \caption{Numerical solution of DDE system for singulary perturbed problem in the case of $\varepsilon = 10^{-2}$.}
    \label{fig:solution-spp-dde-1}
\end{figure}

\begin{figure}
 \centering
    \subfigure[solution~$y_1$]{
    \includegraphics[scale=0.48]{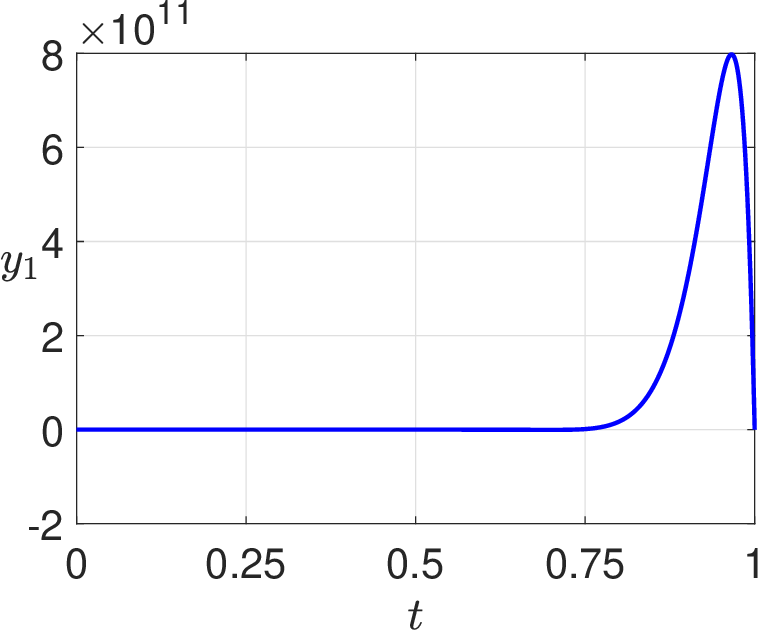}}
    \quad
    \subfigure[solution~$y_2$]{
    \includegraphics[scale=0.48]{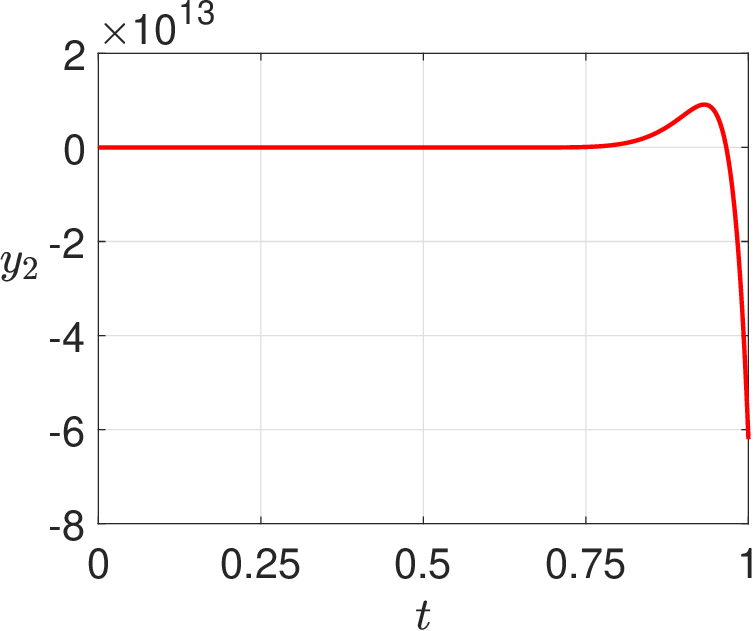}}
    \caption{Numerical solution of DDE system for singulary perturbed problem using $\varepsilon = 10^{-6}$.}
    \label{fig:solution-spp-dde-2}
\end{figure}

Figure~\ref{fig:solution-spp-dde-1} and Figure~\ref{fig:solution-spp-dde-2} 
depict the numerical solutions for $y_1$ as well as $y_2$ 
given two instances of the parameter~$\varepsilon$. 
In the ODE case~\eqref{spp-ode}, the derivatives $u'$ increase rapidly 
for $\varepsilon \rightarrow 0$ at the left boundary ($t=0$).
In contrast, the derivative $u'=y_2$ becomes huge for 
$\varepsilon \rightarrow 0$ at the right boundary ($t=1$) 
in the DDE case~\eqref{spp-dde}. 
Table~\ref{tab:spp-dde} shows the computed boundary values 
(rounded to several digits) and the number of integration steps 
to solve one IVP using the final approximation of the initial value. 
We observe that the singularly perturbed problem was successfully solved using the proposed approach. 
However, the numerical solution process becomes more critical for decreasing values of the parameter~$\varepsilon$.

%%% Figure: SPP - DDE solution %%%%%%%%%%%%%%%%%%%%%%%%%%%%%%%%%%%%%%%%%%%%%%
%%%%%%%%%%%%%%%%%%%%%%%%%%%%%%%%%%%%%%%%%%%%%%%%%%%%%%%%%%%%%%%%%%%%%%%%%%%%%

%%% Figure: SPP - DDE solution %%%%%%%%%%%%%%%%%%%%%%%%%%%%%%%%%%%%%%%%%%%%%%
%%%%%%%%%%%%%%%%%%%%%%%%%%%%%%%%%%%%%%%%%%%%%%%%%%%%%%%%%%%%%%%%%%%%%%%%%%%%%

%%%%%%%%%%%%%%%%%%%%%%%%%%%%%%%%%%%%%%%%%%%%%%%%%%%%%%%%%%%%%%%%%%%%%%%%%%%%%
\begin{table} 
\caption{Boundary values of numerical solutions and number of integration steps in singulary perturbed problem of DDEs.\label{tab:spp-dde}}
\begin{center}
    \begin{tabular}{crccc}
    $\varepsilon$ & $y_2(0)$ & $y_2(1)$ & $y_1(1)$ & no. steps \\ \hline
    $10^{-1}$ & 5.54 & -$2.0 \cdot 10^0$ & 1.0000 & 124 \\
    $10^{-2}$ & -132.00 & -$2.6 \cdot 10^1$ & 1.0000 & 160 \\
    $10^{-3}$ & -38.13 & $4.3 \cdot 10^3$ & 1.0000 & 201 \\
    $10^{-4}$ & -21.37 & $1.5 \cdot 10^6$ & 1.0000 & 259 \\
    $10^{-5}$ & -24.11 & -$2.4 \cdot 10^9$ & 1.0040 & 297 \\
    $10^{-6}$ & -23.60 & -$6.2 \cdot 10^{13}$ & 0.9936 & 379
    \end{tabular}
\end{center}
\end{table}
%%%%%%%%%%%%%%%%%%%%%%%%%%%%%%%%%%%%%%%%%%%%%%%%%%%%%%%%%%%%%%%%%%%%%%%%%%%%%

\clearpage

%%%%%%%%%%%%%%%%%%%%%%%%%%%%%%%%%%%%%%%%%%%%%%%%%%%%%%%%%%%%%%%%%%%%%%%%%%%%%
%%%                          Conclusions                                  %%%
%%%%%%%%%%%%%%%%%%%%%%%%%%%%%%%%%%%%%%%%%%%%%%%%%%%%%%%%%%%%%%%%%%%%%%%%%%%%%

\section{Conclusions and outlook}
We considered a DDE with distributed delay, 
where the weight function of the distribution is a polynomial. 
We derived an equivalent system of DDEs including two discrete delays. 
An exact solution of this system yields an exact solution of the original DDE. 
% provided that the functions are sufficiently smooth. 
We applied a scaling to the DDE with distributed delay, which keeps the magnitude of additional functions small in the equivalent system. 
Alternatively, we arranged an approximating DDE with multiple discrete delays, 
where the integral of the distribution is discretised by Gaussian quadrature. 
IVPs of DDEs with discrete delays were solved by a numerical integrator. 
Results of numerical computations demonstrate a very good agreement of the calculated solutions from the equivalent system and the Gaussian quadrature approach, although an initial discontinuity of the first-order derivative is transferred to discontinuities of higher-order derivatives at later times.  
Thus, solving the equivalent system features three advantages in comparison to an application of a quadrature rule: 
\begin{enumerate}
\item 
A quadrature error is omitted, 
which emerges in the discretisation of the integral that defines the distributed delay. 
\item 
The number of discontinuities in higher-order derivatives of the solution is lower, because a system from a quadrature rule typically involves a larger number of discrete delays. 
\item 
The computation work is often lower, since just two discrete delays appear in the equivalent system, whereas a quadrature rule may include a large number of nodes and thus generates many discrete delays. 
\end{enumerate}
However, the equivalent system includes a weakly unstable part. 
Consequently, errors may increase with time in the numerical solution of IVPs.  
In this case, 
the time step size has to be kept sufficiently small in the 
integration of the equivalent system. 
Yet, the unstable part is not present in the case of DDEs with a constant weight function. 

We outline two possibilities for future work in this context. 
In~\cite{elhameed1,elhameed2}, 
fractional (partial) DDEs were considered, i.e., 
the left-hand side of the differential equations consists of a fractional derivative applied to the state variables. 
If a distributed delay appears on the right-hand side, 
then the approach presented in this article can be investigated to solve such a model. 
A direct application of the approach from Section~\ref{sec:system} would yield a system that includes a non-integer-order derivative of the state variables and a first-order derivative of the auxiliary variables.  
%Furthermore, the construction of equivalent systems can be examined in the case of distributed delays, which are not determined by polynomial weight functions. 
Furthermore, DDEs can be examined, where the right-hand side includes a polynomially distributed delay as well as another distributed or discrete delay. 
A polynomially distributed delay together with a gamma-distributed delay would generate an equivalent system consisting of both DDEs with two discrete delays and ODEs due to the linear chain technique.

%\bigskip\bigskip

%{\bf \Large Statements and Declarations}

%{\bf Funding:} \\ 
%The author states that no funding was received for the work presented in this article.

%{\bf Disclosure of interest:} \\
%The author reports there are no competing interests to declare.

%{\bf Data availability statement:} \\ 
%The data that support the findings of this study are available from the corresponding author, RP, upon reasonable request.

%%%%%%%%%%%%%%%%%%%%%%%%%%%%%%%%%%%%%%%%%%%%%%%%%%%%%%%%%%%%%%%%%%%%%%%%%%%%%
%%%                           References                                  %%%
%%%%%%%%%%%%%%%%%%%%%%%%%%%%%%%%%%%%%%%%%%%%%%%%%%%%%%%%%%%%%%%%%%%%%%%%%%%%%

\end{document}